\documentclass[11pt]{article}
\usepackage{amsmath}
\usepackage{amssymb}
\newcommand{\betti}{\beta}
\newcommand{\comp}{\mathbb{C}}

\newcommand{\id}{{\rm Id}}
\newcommand{\ind}{{\rm ind}}
\newcommand{\kernel}{{\rm ker}}
\newcommand{\La}{\Lambda}
\newcommand{\la}{\lambda}
\newcommand{\mul}{{\rm mul}}
\newcommand{\nulli}{{\rm null}}
\newcommand{\ob}{\overline{b}}
\newcommand{\oc}{\overline{C}}
\newcommand{\q}{\mathbb{Q}}
\newcommand{\qed}{\hfill \Box}
\newcommand{\R}{\mathbb{R}}
\newcommand{\z}{\mathbb{Z}}
\newtheorem{thm}{Theorem}
\newtheorem{pro}{Proposition}
\newtheorem{cor}{Corollary}
\newtheorem{lem}{Lemma}
\newtheorem{rem}{Remark}

\author{Hans-Bert Rademacher}
\title{Existence of closed geodesics on positively 
curved Finsler manifolds} 
\begin{document}
\maketitle
\begin{abstract}
For non-reversible Finsler metrics of positive flag curvature
on spheres and projective spaces
we present results about the number and the length of closed 
geodesics and about their stability properties.
\\ 
\medskip
\\
2000 MSC classification: {\sc 53C22; 53C60; 58E10} 
\end{abstract}
\section{Introduction}
\label{sec:introduction}
For a Finsler metric $F$ on a compact manifold we introduced
in \cite{Ra2004} the concept of {\em reversibility}
$\lambda:=\max \{F(-X)\,;\,F(X)=1\}\ge 1.$
The reversibility attains its minimal value one if and
only if the Finsler metric is {\em reversible,} i.e. 
$F(X)=F(-X)$ for all tangent vectors $X.$ 
In particular Riemannian metrics are reversible.
In this paper we investigate the consequences of the
following length estimate for closed geodesics on
a compact manifold with non-reversible Finsler metric
of positive flag curvature:
\begin{thm}{\rm \cite[Thm.1, Thm.4]{Ra2004}}
\label{thm:length-closed-geodesics}
Let $M$ be a compact and simply-connected manifold with
a Finsler metric $F$ with reversibility $\lambda$ 
and flag curvature $K$ satisfying
$0<K\le 1$ resp. $\frac{\lambda^2}{(\lambda+1)^2}<K\le1$
if the dimension $n$ is odd. Then the length of a closed geodesic
is bounded from below by $\pi \frac{\lambda+1}{\lambda}.$
\end{thm}
It is a generalization of Klingenberg's injectivity radius
estimate for compact Riemannian manifolds, cf. \cite[ch.2.6]{Kl2}.
As announced in \cite[Rem.3]{Ra2004} we apply 
Theorem~\ref{thm:length-closed-geodesics}
to obtain existence results for closed geodesics
on positively curved manifolds carrying a non-reversible
Finsler metric. 
\\
At first we consider the case of the $2$-sphere. In the Riemannian
case there are three geometrically distinct and simple closed geodesics 
on the $2$-sphere  with
length in the interval $[2\pi,2\pi/\sqrt{\delta}]$ if the Gaussian
curvature $K$ satisfies $1/4 < \delta \le K \le 1.$
A closed curve is called {\em simple,} if it does not have
self-intersections.
This is a particular case of an existence result for spheres in all
dimensions, cf. \cite[Thm.A]{BTZ83}.
\\
The {\em Katok metrics} on the $2$-sphere define a one-parameter
family $F_{\epsilon}, \epsilon \in [0,1)$ of 
Finsler metrics of constant
flag curvature $1$ and reversibility $\lambda=(1+\epsilon)/(1-\epsilon),$
for $\epsilon=0$ this is the standard Riemannian metric.
For irrational $\epsilon$ there are exactly two
geometrically distinct closed geodesics $c_1,c_2$ with lengths
$L(c_1)=\pi (1+\lambda^{-1}); L(c_2)=\pi (\lambda +1)\,,$
these geodesics only differ by orientation,
cf. \cite[ch.5]{Ra2004}, \cite[p.142]{Zi}.
\\
Using the Morse inequalities, the topology of the space of
unparametrized closed curves on $S^2$ and a detailed analysis of
the sequence $\ind (c^m)$ of Morse indices of the $m$-fold
covers $c^m$ of a closed geodesic $c$
we obtain the following existence result together with a length
estimate:
\begin{thm}
\label{thm:twosphere}
Let $F$ be a 
Finsler metric on the $2$-sphere with reversibility 
$\la$ and flag curvature $K$ satisfying
$$ \left(\frac{\la}{\la +1}\right)^2< \delta \le K \le 1$$
for some $\delta \in \R^+\,.$ 
Then there are
at least {\em two} geometrically distinct closed geodesics $c_1, c_2$
whose lengths $L(c_1)\le L(c_2)$ satisfy:
$$ L(c_1) \le \frac{2\,\pi}{\sqrt{\delta}}\,;\, L(c_2) \le  
\frac{\pi}{\sqrt{\delta}}\,
\left(\frac{1}{\sqrt{\delta}\frac{\la+1}{\la}-1}+3\right)
\,.$$ 
In addition the shorter closed geodesic $c_1$ is simple.
\end{thm}
\begin{rem} \rm
\label{rem:theorem}
\begin{itemize}
\item[(a)]
If we choose in  particular 
$\delta:=\left(\frac{2 \la +1}{2\la+2}\right)^2
=\left(1-\frac{1}{2\lambda+2}\right)^2$
then \\
$L(c_2) \le 2 \pi (\lambda +2)$
\item[(b)]
If the metric has constant flag curvature, i.e.
$\delta =1$ then $L(c_2)\le \pi (\lambda +3).$
In the above mentioned Katok examples $
L(c_2)=\pi (\lambda+1)/\lambda\,;\,L(c_2)=\pi (\lambda+1)\,.$
A Finsler metric is called {\em bumpy,} if
all closed geodesics are non-degenerate.
Equation~\ref{eq:bdm} implies that for a bumpy Finsler metric
of constant flag curvature $1$ with only two geometrically
distinct closed geodesics $c_1,c_2$ the following relation
holds: $\frac{1}{L(c_1)}+\frac{1}{L(c_2)}=\frac{1}{\pi}\,.$
\item[(c)] The arguments in the Proof of Theorem~\ref{thm:twosphere}
also show: If there is only one 
geometric closed geodesic on the $2$-sphere, then its 
average index is at most $1,$ cf. Theorem~\ref{thm:av}.
It was remarked by Ziller, that a single 
non-degenerate closed geodesic
can generate the homology of the free loop space
up to a any fixed dimension, cf. \cite[p.149]{Zi}.
He also notes that a single closed geodesic $c$ 
with $\ind (c)=1, \ind(c^{2k})=\ind(c^{2k+1})=2k+1, k\ge1$
which 
would produce the homology and for which the energy functional would be
a perfect Morse function, is necessarily degenerate. 
The results of \cite{Ra1} show
that a non-reversible and bumpy
Finsler metric on the $2$-sphere carries at least two geometrically distinct
closed geodesic. More precisely: A bumpy Finsler metric on
the $2$-sphere with only finitely many geometrically 
distinct closed geodesics has
at least two geometrically distinct {\em elliptic} closed geodesics,
cf. \cite[Example 4.1]{Ra1}.
Bangert and Long announced a proof
that for {\em every} non-reversible Finsler metric on the $2$-sphere
there are two geometrically distinct closed geodesics, cf. \cite{Lo1}.
\end{itemize}
\end{rem}
We also present applications of 
Theorem~\ref{thm:length-closed-geodesics} in 
higher dimensions. We use an existence result
for closed geodesics which the author derived in \cite{Ra5}
using the concept of the {\em Fadell-Rabinowitz index}.
We obtain a chain of subordinate cohomology classes in 
a quotient space of the space of closed curves. In the case
of positive flag curvature we can estimate the number of
geometrically distinct closed geodesics whose multiples are
represented by the cohomology classes in this chain.
As a general result we obtain Theorem~\ref{thm:n}
for metrics on compact and simply-connected
manifolds of the rational homotopy type of
a compact rank one symmetric space. Consequences for Finsler
metrics on spheres are listed in the following 
\begin{thm}
\label{thm:nn}
Let $F$ be a 
Finsler metric on the $n$-sphere $S^n$ with reversibility 
$\la$ and flag curvature $K$ satisfying
$0 < \delta \le K \le 1$
for some $\delta \in \R^+\,.$ 
\begin{itemize}
\item[(a)] The number of geometrically distinct
closed geodesics with length $< 2 n\pi$ is
at least $n/2 -1\,,$ 
provided $\delta > \frac{\la^2}{(\la +1)^2}\,.$
\item[(b)] If $\lambda <\frac{n-1}{n-3}, n\ge 4$ and
$\sqrt{\delta}>2\frac{n-2}{n-1}\frac{\lambda}{1+\lambda}$ then
there are at least $(n-2)$ geometrically distinct closed geodesics.
\item[(c)] If $n \ge 6$ is even and $\delta >4/(n-2)^2$ 
then there are at least two geometrically
distinct closed geodesics with length $\le n \,\pi  .$
\end{itemize}
\end{thm}
There are Katok-metrics on the $2k$-sphere $S^{2k}$ resp.
the $(2k-1)$-sphere $S^{2k-1}$ with $2k $ geometrically
distinct closed geodesics (cf. \cite[p.139]{Zi})
and it is an open question whether there are always at least
$n$ closed geodesics on the $n$-sphere \cite[p.156]{Zi}. 
\\[1ex]
Finally we improve this result in Theorem~\ref{thm:bdm}
in the particular case of
a {\em bumpy} metric, i.e. a metric all of whose closed geodesics
are non-degenerate. For the $m$-dimensional complex projective space
$\comp P^m$ we obtain in Corollary~\ref{cor:cpm} 
a lower bound for the number of geometrically
distinct closed geodesics as well as for the 
number of non-hyperbolic closed geodesics provided
there are only finitely many geometrically distinct ones.
\\
One can study stability properties of a closed geodesic $c$ with the
help of the linearized Poincar\'e mapping $P_c\,,$ which is a linear
symplectic map of an $(2n-2)$-dimensional vector space.
It can be defined using the Jacobi fields along this
geodesic, cf. \cite[ch.1]{BTZ}, \cite{Ra2}.
In the most unstable case no (complex) eigenvalue of $P_c$ lies 
on the unit circle, then the closed geodesic is called {\em hyperbolic.}
For example 
all closed geodesics on a Finsler manifold with negative flag curvature
are hyperbolic. 
We obtain a result similar to \cite[Thm.B]{BTZ}:
\begin{thm}
\label{thm:hyperbolic}
Let $F$ be a Finsler metric on a compact manifold
with reversibility $\la$ and flag curvature
$0<\delta\le K\le 1\,.$ There exists a non-hyperbolic closed
geodesic if the $l$--th 
homotopy group $\pi_l(M)$ is non-trivial 
for some $l\ge 2$ and
$\sqrt{\delta} > \frac{l-1}{n-1}\,\frac{\lambda}{\lambda+1}\,.$ 
\end{thm}
In particular on the $n$-sphere $S^n$ with a Finsler metric
satisfying $\lambda^2/(\lambda+1)^2<K\le 1$ there exists a 
non-hyperbolic closed geodesic.
A more detailed analysis also produces existence
results for closed geodesics of elliptic-parabolic type. Here
a closed geodesic is called of {\em elliptic-parabolic type}
if the linearized Poincar\'e map splits into two-dimensional rotations and
a part whose eigenvalues are $\pm 1.$ 
Following the ideas of Thorbergsson \cite{Th} and Ballmann, 
Thorbergsson and Ziller \cite{BTZ}
we obtain as another consequence of the length estimate 
Theorem~\ref{thm:length-closed-geodesics}:
\begin{thm}
\label{thm:elliptic}
On a compact Finsler manifold $M$ with reversibility
$\lambda$ and flag curvature $0<\delta \le K\le 1$
there exists a closed geodesic of elliptic-parabolic type
if one of the following conditions is satisfied:
\begin{itemize}
\item[(a)] $M=S^n$ and $\delta > 
\frac{9}{4} \frac{\la^2}{(\la+1)^2}$ 
with $\la <2\,.$ 
\item[(b)] $M=\R P^n$ and $\delta > \frac{\la^2}{(\la+1)^2}\,.$
\end{itemize}
\end{thm}
It is mentioned in \cite[p.61]{BTZ3} that 
most of the results presented in the Riemannian case
generalize to Finsler metrics. For example it is stated that a
Finsler metric with $9/16<K\le1$ carries a short
closed geodesic of elliptic-parabolic type. But the arguments only work
for {\em reversible} Finsler metrics respectively under the
additional assumption that a shortest closed geodesic
has length $\ge 2\pi\,.$ 
\\
Another setting in which one can show the existence of a closed
geodesic of elliptic-parabolic type is in the presence of
an isometric $S^1$-action.
For example the above mentioned
Katok metrics $F_{\epsilon}$ on the $2$-sphere carry an isometric
$S^1$-action. For irrational parameter $\epsilon$ there are 
exactly two geometrically
distinct closed geodesics which both are elliptic and invariant
under the $S^1$-action. As an analogous result to \cite[Theorem A(iii)]{BTZ}
we obtain:
\begin{thm}
\label{thm:killing}
On a compact manifold with Finsler metric 
with an isometric $S^1$-action 
there exist at least two
geometrically distinct closed geodesics. These closed
geodesics are $S^1$-invariant and they are of elliptic-parabolic
type.
\end{thm}
\section{Critical Point Theory for Closed Geodesics}
\label{sec:critical}
Here we list a couple of results of the critical point theory for 
closed geodesics, general references are the survey article
\cite{Ba2} by Bangert, the book \cite{Kl1} by Klingenberg
and \cite{Ra2}.\\
If $c: S^1:=[0,1]/\{0,1\} \rightarrow M$ is a closed geodesic on the Finsler
manifold $(M,F)$ of length $L(c)$ then for every positive integer
$m$ the $m$-fold cover $c^m: S^1 \rightarrow M; c^m(t)=c(mt)$ is a closed
geodesic, too. If 
$L(c)=\int_0^1 F(c'(t)) \,dt$ denotes the {\em length,} then we have
$L(c^m)=mL(c).$ We call a closed geodesic {\em prime} if it is not
the cover $c_0^m$ of another closed geodesic $c_0$ with $m >1.$ 
Closed geodesics
are the critical points of the energy functional
$$E: \La M \rightarrow \R\,;\, 
E(c)=\frac{1}{2} \int_0^1 F^2\left(c'(t)\right) dt$$
on the {\em Hilbert manifold} $\La M$ of closed curves 
which is the set of all
absolutely continuous closed curves with a square-integrable derivative.
\\
The Morse index $\ind (c)$ of a closed geodesic is the index of the
hessian $d^2E(c)$ of the energy functional. On the space $\La M$ there is a
$S^1$-action given by changing the initial point. The energy functional
is invariant under this group action. 
We call two closed geodesics $c_1,c_2$ of a non-reversible Finsler metric
{\em geometrically equivalent,} if their traces $c_1(S^1) = c_2(S^1)$ and 
their orientation coincide. 
Otherwise we call them {\em geometrically distinct.} In contrast to the
reversible case resp. the case of a Riemannian metric for a closed
geodesic $c$ the curve $c^{-1}$ with $ c^{-1}(t)=c(1-t)$ 
defined by reversing the orientation
in general is not a geodesic. 
Hence a prime closed geodesic $c$ produces
infinitely many critical orbits $S^1.c^m; m \ge 1$ of the energy functional
consisting of all geometrically equivalent
closed geodesics. If a closed geodesic $c$ is the $m$-fold cover 
$c=c_0^m$ of a
prime closed geodesic $c_0$ then we call $m=\mul(c)$ the
{\em multiplicity} of the closed geodesic $c.$ 
Therefore a prime closed geodesic
$c$ produces a {\em tower} $S^1.c^m; m\ge1$ of closed geodesics resp.
critical orbits of the energy functional.
\\
We can view the hessian of the energy functional also as a self-adjoint 
endomorphism. Then the index is the sum of the dimensions of negative
eigenvalues and we call the {\em nullity} $\nulli (c)$ the dimension of
the kernel $\ker d^2E(c)$ minus $1.$ Note that the dimension of the kernel
is always at least $1$ provided $L(c)>0$
since there is a $1$-dimensional group leaving the 
energy functional invariant. A closed geodesic $c$ is called 
{\em non-degenerate} if $\nulli(c)=0.$ Geometrically the nullity
is the dimension of periodic Jacobi fields along $c$ which are orthogonal
to the velocity field $c'.$ Therefore $\nulli (c) \le 2n -2.$
\\
The sequence $\ind (c^m)$ grows almost linearly, we call
the limit
$$\alpha_c:=\lim_{m\to \infty}\frac{\ind (c^m)}{m}$$
introduced by Bott \cite[Cor.1]{Bo}
the {\em average index} and $\overline{\alpha}_c=\alpha_c / L(c)$ the
{\em mean average index.} We have the following estimate for the 
sequence $\ind (c^m):$
\begin{equation}
\label{eq:indcm}
\left| \ind(c^m) - m \alpha_c\right| \le n-1 \,,
\end{equation}
cf. \cite[(1.4)]{Ra1}.
By a Rauch comparison argument as in the Riemannian case one obtains
\begin{lem}
{\rm (\cite[Lem.3]{Ra2004})}
\label{lem:meanaverageindex}
Let $c$ be a closed geodesic on a Finsler manifold $(M,F)$
of dimension $n$ with positive flag curvature $K \ge \delta$
for some $\delta \in \R^+.$
\begin{itemize}
\item[(a)] The mean average index is bounded
from below:
$\overline{\alpha}_c\ge \sqrt{\delta}\, (n-1)/\pi\,.$
\item[(b)] If the length $L(c)$ satisfies 
$L(c) > k  \pi / \sqrt{\delta}$
for some positive integer $k$ then $\ind (c) \ge k (n-1).$
\end{itemize}
\end{lem}
Combining Lemma~\ref{lem:meanaverageindex} with the length estimate
Theorem~\ref{thm:length-closed-geodesics} 
for a closed geodesic we obtain:
\begin{lem}
\label{lem:averageindex}
Let $c$ be a closed geodesic on a compact and
simply-connected Riemannian manifold of dimension $n$
with a non-reversible
Finsler metric with reversibility $\lambda$ and flag curvature
$0 < \delta \le K \le 1$ 
where $\delta > \frac{\lambda^2}{(\lambda+1)^2} $
if $n$ is odd. Then
$$\alpha_c \ge \sqrt{\delta} \,\frac{\lambda+1}{\lambda} \,\left(n-1\right).$$
\end{lem}
{\sc Proof.} Since $L(c) \ge \pi \left(1 + 1/\lambda\right)$ 
by Theorem~\ref{thm:length-closed-geodesics} the claim follows
from Lemma~\ref{lem:meanaverageindex}. $\qed$
\\
Now we come to the {\em Morse Inequalities} of the $S^1$-invariant
functional $E : \La M \rightarrow \R.$ Let 
$$\ob_j:= b_j\left(\La M/S^1, \La^0M/S^1; \q\right)$$
where for $a\ge0$ we denote
$\La^aM:=\left\{\sigma \in \La M\,|\, E(a) \le a\right\}$
the sublevel sets and
$b_j$ is the $j$-th Betti number. 
In particular $\La^0M$ is the set of point curves
which can be identified with the manifold $M.$ Since it is the
fixed point set of the $S^1$-action one can also identify
the quotient space 
$\La^0 M/S^1$ with the manifold $M.$ Given a closed geodesic
$c$ we use the following notation
$$\La (c):= \left\{ \sigma \in \La M \,|\, E(\sigma)< E(c)\right\}.$$
Then we call
$$\oc_*(c)= H_*\left(\left(\La (c) \cup S^1.c\right)/S^1,\La (c)/S^1; \q\right)\,$$
the $S^1${\em -critical group} of the closed geodesic $c$
and let $\ob_j(c)=\dim \oc_j(c).$
\\
We collect the information about the $S^1$-critical groups in the following 
two lemmas:
\begin{lem}
\label{lem:critical-group-nondegenerate}
{\rm (\cite[Proof of Prop.2.2]{Ra1})}
Let $c$ be a {\em non-degenerate} closed geodesic with
$i=\ind (c), m=\mul(c)\,.$ Then $c=c_1^m$ for a prime closed
geodesic $c_1$ and 
$$
\ob_j(c)= \left\{\begin{array}{lll}
1 & ; & j=i; m\equiv 1 \pmod{2}\\
1 & ; &  j=i; m\equiv 0 \pmod{2}\\
&& \mbox{ \rm and }\ind(c_1^2) \equiv \ind (c_1) \pmod{2} \\
0 &;& \mbox{ \rm otherwise }
\end{array}\right.
$$
\end{lem}
In the general case we obtain the following cases:
\begin{lem}
\label{lem:critical-group-general-case}
{\rm \cite[Satz 6.13]{Ra2}}
Let $c$ be a closed geodesic with $i=\ind (c); l=\nulli(c).$ 
Then we have the following statements:
\begin{itemize}
\item[(a)] $\ob_j(c)=0$ for $j<i$ or $j>i+l.$
\item[(b)] $\ob_i(c) +\ob_{i+l}\le 1$ and if $\ob_i(c)+\ob_{i+l}=1$ then
$\ob_j(c)=0$ for all $j$ with $i+1 \le j\le i+l-1.$
\end{itemize}
\end{lem}
As a consequence of the formula for the sequence $\ind (c^m), m \ge 1$
given by Bott \cite{Bo} we conclude:
\begin{lem}
\label{lem:parity-indcm}
Let $c$ be a closed geodesic on a surface (i.e. $n=\dim M=2$) 
with $\ind(c)=1$
and average index $\alpha_c >1.$ Then for all
$m\ge 1$ the indices $\ind (c^m)$ are odd.
\end{lem}
{\sc Proof.}
Let $P_c$ be the linearized Poincar\'e mapping, i.e. the 
linearization
of the return map of the closed orbit of the geodesic flow
corresponding to the closed geodesic. 
There is a function
$$I_c: S^1:=\left\{z \in \comp \,;\, |z|=1\right\} 
\rightarrow \z^{\ge 0}$$
with the following properties, 
cf. \cite[Thm. A,C]{Bo},\cite[ch.9]{Lo},\cite[ch.4]{Ra2}:
\begin{itemize}
\item[(a)]
$\ind (c^m) = \sum_{z^m=1} I_c(z) $
\item[(b)] Define $N_c: S^1 \rightarrow \z^{\ge 0}:$
$N_c(z) = \dim \kernel (P_c - z \id)\,.$ Then
$$\nulli (c^m) = \sum_{z^m=1} N_c(z)\,.$$
\item[(c)]
The function
$I_c$ is constant in a neighborhood of points $z$ with
$N_c(z)=0.$ For the {\em splitting numbers}
$$
S_c^{\pm}(z) ;= \lim _{\phi \to \pm 0} I_c(z \exp(i \phi)) - I_c(z)
$$
of the function $I_c$ the following estimate holds:
$$
0 \le S_c^{\pm }(z) \le N_c (z)
$$
\item[(d)]
$
I_c(z)=I_c(\bar z), N_c(z) = N_c(\bar z)
$
\end{itemize}
It follows that 
$\ind(c^m) \equiv \ind(c)\pmod{2}$ for all odd $m$ and
$\ind(c^m) \equiv I_c(1)+I_c(-1)=\ind(c^2)\pmod{2}$ for all even
$m.$ It was also shown by Bott that the splitting numbers
only depend on the symplectic normal form of the
linearized Poincar\'e map $P_c,$ for a detailed discussion
see \cite[(2.13)]{BTZ}, \cite[ch.IV]{Lo}, \cite[ch.4]{Ra2}.\\
Now we come to the case $n=2,$ then for the eigenvalues $z$
of the linearized Poincar\'e map there are the 
following cases:
\begin{itemize}
\item[(a)]
$z \not\in S^1,$ i.e. $z$ is a real number
with $|z|\not=1.$ Then also $z^{-1}$ is an eigenvalue,
in this case the closed geodesic is called 
{\em hyperbolic} and $\ind (c^m)=m \,\ind(c).$
In particular the average index satisfies 
$\alpha_c=1$ in contradiction to our assumption.
\item[(b)]
If $z=1$ then $S^+(1)=S^-(1)$ and
$\alpha_c=I_c(-1)=\ind (c^2)-\ind (c)=\ind(c)+S^+(1)=1+S^+(1)>1$ 
by assumption, hence
we conclude $S^+(1)=1.$ It follows that $\ind(c^2)=3$
and $\alpha_c=2.$
\item[(c)]
If $z=-1$ then 
$\alpha_c=\ind(c)=1$ in contradiction to our assumption.
\item[(d)]
If $z=\exp(\sqrt{-1}\pi \rho)$ with $\rho \in (0,1)$ then
we conclude from \cite[(2.13)]{BTZ} or \cite[Thm.4.3]{Ra2}:
$S^{+}(z)+S^{-}(z)=1.$ Since $\alpha_c=\ind(c) + (S^+(z) - S^-(z))\rho >1$
we conclude $S^+(z)=1, S^-(z)=0, \alpha_c\in(1,2)$ 
Hence in this case $\ind(c^2)=I_c(1)+I_c(-1)=
2\,\ind(c)+S^+(z)=3.$
\end{itemize} 
Therefore $\ind (c^2)=3$ which implies that for all even $m$
the indices $\ind(c^m)$ are odd, too. 
$\qed$
\\
Now Lemma~\ref{lem:critical-group-general-case} and
Lemma~\ref{lem:parity-indcm} imply the following
\begin{cor}
\label{cor:dimension-two}
Let $c$ be a prime
closed geodesic on a surface (i.e. $n=2$) with
Finsler metric with index $\ind(c)=1$ and
average index $\alpha_c>1\,.$ 
Then for every $m\ge1:$
$$\sum_{j \equiv 1\pmod{2}} \ob_{j}(c^m) \le 1\,.$$
\end{cor}
The Morse inequalities relate the critical groups as local information
about the critical points with the global topological information given
by the Betti numbers of the space on which the Morse function is defined.
\begin{lem}
\label{lem:betti-two-sphere}
{\rm (\cite[2.6]{Ra1})}
The rational Betti numbers 
$\betti_i:=b_i\left(\La S^2/S^1,\La^0 S^2/S^1; \q\right),$ 
of the pair of
quotient spaces $\left(\La S^2/S^1,\La^0 S^2/S^1\right)$ are given by:
$$ \betti_i = \left\{
\begin{array}{ccl}
2 &;&i=2m+1, m \ge 1\\
1&;& i=1\\
0 &;& i=2m, m\ge 0
\end{array}
\right.
$$
\end{lem}
\section{Proof of Theorem~\ref{thm:twosphere}}
Let $N$ be the odd integer satisfying
\begin{equation}
\label{eq:N} 
N-2 \le 
\frac{1}{\sqrt{\delta}\frac{\la+1}{\la}-1} < N\,.
\end{equation}
We assume that there is only one 
class of geometrically equivalent closed geodesics whose indices
are bounded from above by $N.$ 
Hence there is a prime closed geodesic $c$ 
such that every closed
geodesic $\tilde{c}$ with $\ind(\tilde{c})\le N$ is up to the 
choice of the initial point of the form
$c^m\,$ respectively $\tilde{c}\in S^1.c^m$ for some 
$m\ge1\,.$
\\
We define for all $i$
with $0 \le i\le N: \,v_i:= \sum_{m\ge 1} \ob_i(c^m),$ 
then the Morse Inequalities
for the $S^1$-invariant energy functional 
$E: \La S^2 \rightarrow \R$
yield (cf. \cite[ch.6.1]{Ra2}):
\begin{equation}
\label{eq:morseinequalities}
v_i \ge \betti_i
\end{equation}
for all $i$ with $0 \le i\le N.$
In particular we conclude from $\betti_1=1$
that for some $m \ge 1:\,
\ind (c^m) \le 1.$ Since 
$L(c) \ge \pi \frac{\la+1}{\la}$ and 
$K\ge \delta > \left(\frac{\la}{\la+1}\right)^2$
we obtain from Lemma~\ref{lem:meanaverageindex} that 
$\ind (c^m) \ge 1$ for all $m \ge 1.$ Hence  
we have finally shown: $\ind (c)=1.$ 
As an estimate for the average index 
we obtain from Lemma~\ref{lem:meanaverageindex} :
$
\alpha_c \ge \sqrt{\delta}\,(\la+1)/\la >1\,.
$
Inequality (~\ref{eq:indcm}) and Inequality (~\ref{eq:N})
imply:
$$
\ind(c^N) \ge N \alpha_c-1\ge N \sqrt{\delta}\,\frac{\la+1}{\la} -1
= N \left(\sqrt{\delta}\,\frac{\la+1}{\la}-1\right) +N-1 > N \,.
$$
Therefore Corollary~\ref{cor:dimension-two} implies:
\begin{equation*}
\sum_{
\begin{array}{c}\scriptstyle 0\le i\le N \\
\scriptstyle i \equiv 1\pmod{2}
\end{array}}
v_i  =
\sum_{
\begin{array}{c} \scriptstyle 0\le i\le N ; m\ge 1 \\
\scriptstyle i \equiv 1\pmod{2}
\end{array}}
 \ob_i (c^m)  
\le \#\left\{\left. m\,\right|\, \ind (c^m)\le N\right\} < N\,.
\end{equation*}
This contradicts the Morse Inequalities (~\ref{eq:morseinequalities})
since by Lemma~\ref{lem:betti-two-sphere} for $N$ odd:
\begin{equation*}
\sum_{
\begin{array}{c}\scriptstyle 0\le i\le N \\
\scriptstyle i \equiv 1\pmod{2}
\end{array}}
v_i  \ge 
\sum_{\begin{array}{c} \scriptstyle 0\le i \le N\\
\scriptstyle i \equiv 1\pmod{2}
\end{array}}
\betti_i= N\,.
\end{equation*}
Hence there are at least
two geometrically distinct closed geodesics 
$c_1,c_2$ with $\ind (c_1)=1\,;\, \ind(c_2)\le N\,.$
We conclude from Lemma~\ref{lem:meanaverageindex}
and Inequality (~\ref{eq:N}):
$$
L(c_2)\le \frac{\pi}{\sqrt{\delta}}(N+1) \le
\frac{\pi}{\sqrt{\delta}}
\left(\frac{1}{\sqrt{\delta}\frac{\lambda+1}{\lambda}-1}+3\right)
$$
Since $\ind (c_1)=1$
Lemma~\ref{lem:meanaverageindex}(b) implies that
$L(c_1) \le 2\pi/\sqrt{\delta}\,.$
For the given curvature bounds not only the length
of a closed geodesic but also the length of a geodesic
loop is bounded from above by $\pi (1+\lambda^{-1}),$
cf. \cite[Thm. 1]{Ra2004}, hence $c_1$ is simple
since $2\pi /\sqrt{\delta}< 2\pi (1+\lambda^{-1}).$
$\qed$\\
The assumption about the flag curvature was only used in the 
proof to obtain a lower bound $>1$ for the average index and to
show that the shortest closed geodesic is simple. Hence the proof
also shows the following 
\begin{thm}
\label{thm:av}
Let $F$ be a non-reversible Finsler metric on $S^2$ with a closed
geodesic with average index $\alpha_1=\alpha_{c_1}>1\,.$ Then there
is a second closed geodesic $c_2$ with index
$\ind(c_2)\le \frac{1}{\alpha_1-1}+2\,.$
\end{thm}

\section{Existence results in higher dimensions}
We consider a compact and simply-connected manifold $M$ whose
rational cohomology algebra is generated by a single element 
$x \in H^d(M;\q)$ of
degree $d,$ with the relation
$x^{m+1}=0\,.$ 
Hence the cohomology algebra $H^*(M;\q)$
is isomorphic to the truncated polynomial algebra
$ \cong T_{d,m+1}(x)= \q [x]/(x^{m+1}=0)$
and the dimension of $M$ equals $d \,m.$
The main examples are the compact rank one
symmetric spaces, i.e. spheres $S^d$ of dimension $d$ (then
$m=1$), $m$-dimensional complex projective spaces 
$\comp P^m$ with
$d=2,$ $m$-dimensional quaternionic projective 
spaces $\mathbb{H}P^m$ with
$d=4$ and the Cayley plane $\mathbb{C}{\rm a}P^2$ with
$d=8, m=2.$ Then we obtain from \cite[Thm.5.11]{Ra5}:
\begin{pro}
\label{pro:fadellrabinowitz}
Let $M$ be a simply-connected and compact manifold
whose rational cohomology algebra is generated by a single
element of order $d,$ i.e. $H^*(M;\q)=T_{d,m+1}(x)$
endowed with a Finsler metric. Then there is a sequence $c_k; k\ge 1$ of
prime closed geodesics and a sequence $m_k;k\ge 1$ of 
positive integers such that the sequence $S^1.c^{m_k}_k; k \ge 1$
is a sequence of $S^1-$orbits of closed geodesics which are
pairwise distinct (although in general not geometrically
distinct) and whose lengths and indices 
satisfy the following proerties for all $k \ge 1\,:$
properties:
\begin{itemize}
\item[(a)] 
$m_k L(c_k)=L\left(c_k^{m_k}\right) \le 
L\left(c_{k+1}^{m_{k+1}}\right)=m_{k+1} L(c_{k+1})$
\item[(b)] 
$2k - (2m-1)d+1 \le \ind\left(c^{m_k}_k\right) \le 2k +d-1$
\end{itemize}
\end{pro}
\begin{thm}
\label{thm:n}
Let $M$ be a simply-connected and compact manifold
whose rational cohomology algebra is generated by a
single element $x \in H^d(M;\q)$, i.e. $H^*(M;\q)=T_{d,m+1}(x)\,.$
We assume that the
manifold $M$ carries a Finsler metric with reversibility
$\lambda$ whose flag curvature $K$ is positive and satisfies
$$ 0 < \delta \le K \le 1\,,$$
where $\sqrt{\delta}> \lambda/(\lambda+1)$ if $n$ is odd.
Then the number of geometrically distinct closed geodesics
of length $\le L$ is bounded from below by
$$ A(m,d,\delta, \lambda, L):= 
\frac{1}{2}\frac{\lambda+1}{\lambda}\,\sqrt{\delta}\,
\left(md-1 - \frac{\pi}{L} d\right)
$$
\end{thm}
\begin{rem}
\rm
\begin{itemize}
\item[(a)] If we are not interested in the length of the
closed geodesics we obtain as bound:
$$A(m,d,\delta,\lambda,\infty)=
\lim_{L \to \infty} A (m,d,\delta,\lambda,L)=
\frac{1}{2}(md -1)\,\sqrt{\delta}\,
\frac{\lambda +1}{\lambda}$$
\item[(b)] The maximal value of the bound is attained if
the flag curvature $K$ is constant ($\delta=1$) and 
the metric is reversible ($\lambda=1$):  
$$A (m,d,1,1,\infty)= md-1 = n-1 $$
\item[(c)] Theorem~\ref{thm:nn} is a direct consequence 
for $m=1, d=n\,.$
\end{itemize}
\end{rem}
{\sc Proof of Theorem~\ref{thm:n}:}
We consider the sequence $\left(S^1.c^{m_k}_k\right)_{k\ge1}$ 
of pairwise distinct
critical orbits of closed geodesics satisfying the properties of
Proposition~\ref{pro:fadellrabinowitz}. 
Hence $c_k,k\ge 1$ are prime closed geodesics and 
$(m_k)$ is a sequence of positive integers. Fix a number $L>0$
and let $a_L:=\#\{k\ge 1\,|\, L(c_k)m_k\le L\}\,.$
By the comparison result part (b) 
of Lemma~\ref{lem:meanaverageindex}:
$$a_L \ge b_L:=\#\left\{k\ge 1\,\left|\, \ind \left(c_k^{m_k}\right)\right.<
\frac{L}{\pi}\,\sqrt{\delta}\,(n-1)\right\}\,.$$
We conclude from Proposition~\ref{pro:fadellrabinowitz}, (b)
that for $2k <L\,\sqrt{\delta}\,(n-1)\pi^{-1} -(d-1)\,:\,
\ind(c_k^{m_k})< L \,\sqrt{\delta}\,(n-1) \pi^{-1}\,,$
hence
$b_L$ is bounded from below by the integer part
of $L\,\sqrt{\delta}\,(n-1)\pi^{-1} -(d-1),$ hence
$b_L \ge  L\,\sqrt{\delta}\,(n-1)\pi^{-1} -d\,.$
Since $L(c_k)\ge \pi (\lambda+1)/\lambda$ by 
Theorem~\ref{thm:length-closed-geodesics}
we obtain from Proposition~\ref{pro:fadellrabinowitz}, (a)
that $m_k < L \pi^{-1} \lambda/(\lambda+1)\,.$
Therefore the number of geometrically distinct closed geodesics
in the set $c_1,c_2,\ldots,c_{a_L}$ is bounded from below
by
\begin{eqnarray*}
\frac{a_L}{\frac{L}{\pi} \frac{\lambda}{\lambda+1}}\ge
\frac{\frac{L}{\pi}\,\sqrt{\delta}\,(n-1) -d}
{2 \,\frac{L}{\pi} \frac{\lambda}{\lambda+1}}
\ge \frac{1}{2}\frac{\lambda+1}{\lambda}\,\sqrt{\delta}\,
\left(n-1 - \frac{\pi}{L}\,d\right)
\end{eqnarray*}
$\qed$\\
We call a Finsler metric {\em bumpy}, if all closed geodesics are
non-degenerate. If the Finsler metric on a compact and
simply-connected manifold is bumpy and has only finitely many 
geometric distinct prime closed
geodesics $c_1, c_2,\ldots,c_r$ with average indices $\alpha_1,
\alpha_2,\ldots, \alpha_r$ then the rational
cohomology ring $H^*(M;\q)$ is generated by a single
element $x$ of degree $d$ with the only relation
$x^{m+1}=0,$ i.e. $n=\dim M=m d\,.$ 
The invariants $d,m$ determine the number
$$
B(d,m) = \left\{
\begin{array}{cl} -\frac{m(m+1)d}{2d(m+1)-4} & ; d \mbox{ even } \\ & \\
\frac{d+1}{2d-2} &; d \mbox{ odd }
\end{array} \right. \, .
$$
for which the following formula is derived in
\cite[Thm.3]{Ra1}:
\begin{equation}
\label{eq:bdm}
B(d,m)=\sum_{i=1}^{r} \frac{\gamma_i}{\alpha_i}
\end{equation}
Here $\gamma_i \in \{\pm 1/2,\pm 1\}$ is an invariant controling
the parity of the sequence $\ind (c_i^m)$ and the orientability of the
negative normal bundle $c$ and $c^2\,.$ 
Let $c_1,\ldots,c_s; s\le r$ be the non-hyperbolic closed geodesics.
Then for even $d$ the following estimate holds with the
same argument as in \cite[Thm.3.1(b)]{Ra1}:
\begin{equation}
\label{eq:eps}
\sum_{k=1}^{s}\left\{|\gamma_k|
\left(\frac{md-1}{\alpha_k}-1\right)+2\right\}
\ge \frac{1}{4} m (m+1) d
\end{equation}
As a consequence from this formula and Lemma~\ref{lem:averageindex}
we obtain analogous to \cite[Cor.3.4]{Ra1}:
\begin{thm}
\label{thm:bdm}
Let $F$ be a Finsler metric on a compact and simply-connected manifold 
$M$ with $H^*(M;\q)=T_{d,m+1}(x)$ with reversibility
$\la$ and flag curvature $K$ satisfying
$0<\delta\le K\le 1,$ where $ \sqrt{\delta} > \la/(\la+1)$ provided
$m=1$ and $d$ is odd. 
\begin{itemize}
\item[(a)]
If the metric is bumpy then there are at least\\
$C(m,d,\delta,\la):=
|B(d,m)| \sqrt{\delta}\frac{\la+1}{\la} (md-1)$ geometrically distinct
closed geodesics.
\item[(b)] 
If the metric is bumpy and there exist only finitely many geometrically
distinct closed geodesics then there are at least
$m (m+1) d\left(\frac{4 \lambda}{\sqrt{\delta}(\lambda+1)}+6\right)^{-1}$
non-hyperbolic closed geodesics.
\end{itemize}
\end{thm}
For $d \in\{2,4,8\}$ and for a fixed value 
of the lower curvature bound
$\delta$ the function $C(m,d,\delta,\la)$ grows quadratically
in $m.$ Let us consider Finsler metrics on the 
$m$-dimensional complex procjective space $\mathbb{C}P^m.$
The flag curvature of the
normalized Fubini-Study metric on $\comp P^m$ satisfies
$1/4 \le K \le 1\,.$ 
For $\delta \rightarrow 1/4$ and
$\lambda \rightarrow 1$ we obtain as maximal value 
$C(m,2,1/4,1)= (m+1)(m-1/2)\,.$ 
It is very likely that this bound is not 
optimal since there are Finsler metrics of Katok type on
$\comp P^m$ with $m(m+1)$ geometrically distinct closed geodesics,
cf. \cite[p.139]{Zi}.\\
Another application is the following result
analogous to \cite[Cor.4]{BTZ}:
\begin{cor}
\label{cor:cpm}
A bumpy Finsler metric on the 
$m$-dimensional complex projective space $\comp P^m$ ($m\ge 7$)
with reversibility $\lambda$ and flag curvature
$0<\delta\le K\le 1; \sqrt{\delta}=
\frac{2}{m+1}\frac{\lambda}{1+\lambda}$ 
with only finitely many geometrically distinct closed geodesiscs
carries at least
$2m$ geometrically distinct closed geodesics.
At least $(m-3)$ of these closed geodesics are non-hyperbolic. 
\end{cor}
\section{Stability properties of closed geodesics}
\label{sec:stable}
A closed geodesic is called {\em hyperbolic} if all
all eigenvalues of the linearized Poincar\'e map have
modulus $\not=1\,.$ Then the sequence $\ind (c^m)$ is linear in $m,$
i.e. $\ind (c^m)=m \ind(c)\,.$ This was observed by Bott in
\cite{Bo} and shows immediately part (a) of the following
\begin{lem}
\label{lem:hyper}
Let $c$ be a closed geodesic of a Finsler manifold with
average index $\alpha_c\,.$
\begin{itemize}
\item[(a)] If $\ind (c) \not= \alpha_c$ then $c$ is non-hyperbolic.
\item[(b)] If $\ind (c^2)-2\,\ind(c)=n-1$ then $c$ is of elliptic-parabolic
type, the linearized Poincar\'e map splits into $2 \times 2$ blocks
of the form
$$ \left(\begin{array}{rr}
1 &0\\1&1
\end{array}\right)
\enspace \text{and / or} \enspace
\left(
\begin{array}{rr}
\cos \phi & -\sin \phi\\
\sin \phi & \cos \phi
\end{array}
\right)
\enspace 0\le \phi < \pi\,,
$$
with respect to a symplectic basis 
$(X_1,Y_1,X_2,Y_2,\ldots,X_{n-1},Y_{n-1})$
satisfying $\omega (X_i,Y_i)=\delta_{ij}; \omega(X_i,X_j)=
\omega(Y_i,Y_j)=0$ for the symplectic form $\omega\,.$
\end{itemize}
\end{lem}
Part (b) is shown in \cite[Lemma 3.1]{BTZ}, since $\ind(c^2)=I_c(1)+I_c(-1)$
and $\ind (c)=I_c(1)\,,$ cf. the Proof of Lemma~\ref{lem:parity-indcm}
\\[1ex]
{\sc Proof of Theorem~\ref{thm:hyperbolic}:}
By considering the universal covering we can assume $M$ to be simply-connected.
Since $\pi_l(M)\not=0$ one concludes that $\pi_{l-1}(\Lambda M)\not=0,$
hence there is a closed geodesic $c$ with $\ind (c)\le l-1\,.$
From Theorem~\ref{thm:length-closed-geodesics} we conclude
$L(c)\ge \pi(1+\lambda^{-1}),$ hence Lemma~\ref{lem:meanaverageindex} (b)
implies that $\alpha_c> l-1\,.$ Therefore $c$ is non-hyperbolic by
Lemma~\ref{lem:hyper}, part (a).
$\qed$\\[1ex]
\begin{rem}
\rm On a compact and not-simply-connected Riemannian manifold 
of non-negative Ricci curvature there is a non-hyperbolic closed
geodesic. The proof of this statement (cf. \cite[Thm.B (a)]{BTZ}) 
carries over to the Finsler case without changes, here one does not
need Klingenberg's injectivity radius estimate.
\end{rem}
{\sc Proof of Theorem~\ref{thm:elliptic}:}
\\
(a)\,
Since $b_{n-1}(\Lambda S^n,\Lambda^0 S^n)=1$ there
is a closed geodesic $c$ with $\ind(c)\le n-1\,.$
The length estimate Theorem~\ref{thm:length-closed-geodesics}
implies for the second cover $L(c^2)\ge 2\pi (1+\lambda^{-1})> 
3 \pi/\sqrt{\delta}.$
Hence by Lemma~\ref{lem:meanaverageindex} 
we obtain $\ind(c^2)\ge 3(n-1)$ resp. $\ind(c^2)-2\,\ind(c)\ge n-1.$
We conclude from Lemma~\ref{lem:hyper}(b) that $c$ is 
of elliptic-parabolic type.
\\[1ex]
(b)\,
Let $c$ be a shortest closed geodesic which is not homotopically
trivial. Then $c$ is a local minimum for the energy functional,
hence $\ind(c)=0\,.$ Since $c$ defines also a closed geodesic on
the universal covering we conclude from 
Theorem~\ref{thm:length-closed-geodesics}
that $L(c^2)\ge \pi(1+\lambda^{-1})>\pi/\sqrt{\delta}.$
Hence by Lemma~\ref{lem:meanaverageindex}(b) we conclude 
$\ind(c^2)\ge(n-1)$ which shows that $c$ is of elliptic-parabolic
type by Lemma~\ref{lem:hyper}(b). 
$\qed$\\
Now we come to a different setting in which one can show
the existence of {\em two} geometrically distinct closed
geodesics on a manifold with non-reversible Finsler metric.
\\
An isometry $A:M\rightarrow M$ of finite order 
of a compact Finsler manifold
$(M,F)$ has {\em small displacement,} if for all points $p \in M$
the image point 
$A(p)$ does not lie in the cut locus. 
Let $\theta: M \times M \rightarrow \R$ be the distance function,
i.e. $\theta(x,y)$ is the minimal length of a 
smooth curve $c:[0,1]\rightarrow M$ from $x=c(0)$ to $y=c(1)\,.$
Note that this distance function in general is not symmetric
since we consider non-reversible Finsler metrics.
Then we define 
the function $f_A: M\rightarrow \R$
with $f_A(x)=\theta^2(x, Ax)$ which is smooth outside
the fixed point set $M_A$ of $A\,.$ Then we obtain as in
\cite[p.239]{BTZ}: The point $p$ is a critical point of
$f_A$ if and only if the unique minimal geodesic
$\gamma:[0,1]\rightarrow M$ from
$p=\gamma(0)$ to $A(p)=\gamma(1)$ is invariant under $A,$ 
i.e. $A_*\left(\gamma'(0)\right)=\gamma'(1)\,.$ 
It also follows that $\gamma(t)$ is a critical point of
$f_A$ for all $t\in [0,1].$
A critical point $p \not\in M_A$
determines a closed geodesic, since $A$ is of finite order.
Analogous to \cite[Thm.3.10]{BTZ} we obtain
\begin{pro}
\label{pro:isom}
Let $A$ be an isometry of finite order and small displacement
on a compact Finsler manifold.
A local maximum of $f_A$ determines a closed geodesic of elliptic-parabolic
type if $2 \theta(p,A(p))<\theta(p,{\rm Cut}(p))\,.$
\end{pro}
{\sc Proof of Theorem~\ref{thm:killing}:}
\\
Let $\phi_t, t\in \R$ define the isometric $S^1$-action, i.e.
$\phi_t: M \rightarrow M$ is a one-parameter group of isometries
with $\phi_1$ equals the identity. Then for a sufficiently large
integer $m$ the isometries $A=\phi_{1/m}$ and $A^{-1}=\phi_{-1/m}$ have small
displacement and are of finite order. Then we conclude from
Proposition~\ref{pro:isom} that there are two points $p_{\pm}$ such
that $p_{\pm}$ is a maximum of $f_A$ resp. $f_{A^{-1}}\,.$
Then the closed geodesics $c_{\pm}$ with $c_{\pm}(0)=p_{\pm}$ and
$c_{\pm}(1/m)=A^{\pm 1}(p_{\pm})$ are of elliptic-parabolic type and
satisfy $c(t)=\phi_t(p_{\pm})\,.$
If $c_-, c_+$ are
geometrically equivalent, we can assume without loss of generality
that $p=p_-=p_+\,.$ Hence $c_{\pm}(t)=\phi_{\pm t}(p)$ which shows
that $c_+'(0)=-c_-'(0)\,.$ i.e. the closed geodesics $c_,c_+$ are
geometrically distinct.
$\qed$\\
One can obtain the existence of several invariant and closed geodesics
under additional assumptions:
\begin{pro}
\label{pro:projective}
On the odd-dimensional sphere $S^n$ with $n=2m+1$ endowed with a 
non-reversible Finsler metric with a free and isometric 
$S^1$-action for which we can identify 
the quotient space $S^{2m+1}/S^1$ with the
complex projective space $\comp P^{m}$ there 
are $n+1=2m+2$ geometrically distinct and $S^1$-invariant closed geodesics.
\end{pro}
The number $2m+2=n+1$ is optimal
in this case as certain Katok examples show, cf. \cite[p. 139]{Zi}.
In the proof we use the function
$F_A: \comp P^{m}\rightarrow \R^+$ resp. 
$F_{A^{-1}}:\comp P^{m}\rightarrow \R^+$
induced by the $S^1$-invariant functions
$f_A: S^n\rightarrow \R^+$ resp.
$f_{A^{-1}}:S^n\rightarrow \R^+\,$ as in the Proof of the
preceding
Proposition~\ref{pro:isom}. 
The critical points of $F_A, F_{A^{-1}}$ correspond to (non-trivial)
invariant closed geodesics on $S^{2m+1}$ and $m+1$ is the
cup length of the complex projective space. The critical points of
$F_A$ and $ F_{A^{-1}}$ are geometrically distinct, since they differ
by orientation if they have the same trace.

{\sc Universit\"at Leipzig, Mathematisches Institut\\ Augustusplatz 10/11,
D-04109 Leipzig}\\
{\tt rademacher@math.uni-leipzig.de}\\
{\tt www.math.uni-leipzig.de/\symbol{126}rademacher}
\end{document}